\documentclass{article}
\usepackage{amssymb}
\usepackage{amsfonts}
\usepackage{amsmath}
\usepackage{graphicx}
\usepackage{anysize}
\usepackage{url}

\newtheorem{theorem}{Theorem}

\newtheorem{definition}[theorem]{Definition}
\newtheorem{example}[theorem]{Example}

\newtheorem{lemma}[theorem]{Lemma}

\newtheorem{proposition}[theorem]{Proposition}

\newenvironment{proof}[1][Proof]{\noindent\textbf{#1.} }{\ \rule{0.5em}{0.5em}}
\newenvironment{remark}[1][Remark]{\noindent\textbf{#1.} }{}

\begin{document}

\title{Borsuk-Ulam Theorems for Complements of Arrangements}
\author{Pavle V.~M.~Blagojevi\'{c}, Aleksandra S.~Dimitrijevi\'{c} Blagojevi%
\'{c} \and John McCleary}
\maketitle

\begin{abstract}
In combinatorial problems it is sometimes possible to define a $G$%
-equivariant mapping from a space $X$ of configurations of a system to a
Euclidean space $\mathbb{R}^{m}$ for which a coincidence of the image of
this mapping with an arrangement $\mathcal{A}$ of linear subspaces insures a
desired set of linear conditions on a configuration. Borsuk-Ulam type
theorems give conditions under which no $G$-equivariant mapping of $X$ to
the complement of the arrangement exist. In this paper, precise conditions
are presented which lead to such theorems through a spectral sequence
argument. We introduce a blow up of an arrangement whose complement has
particularly nice cohomology making such arguments possible. Examples are
presented that show that these conditions are best possible.
\end{abstract}


%




\section{Borsuk-Ulam type results}

\label{intro}

\noindent Theorems of Borsuk-Ulam type present conditions preventing the
existence of certain equivariant mappings between spaces. The classical
Borsuk-Ulam theorem, for example, treats mappings of the form $f\colon
S^{n}\rightarrow \mathbb{R}^{n}$ for which $f(-x)=-f(x)$, that is, $f$ is
equivariant with respect to the antipodal action of $\mathbb{Z}/2$ on $S^{n}$
and the action of $\mathbb{Z}/2$ on $\mathbb{R}^{n}$. Such a map \emph{must}
meet the origin.

\noindent Generalizations of the Borsuk-Ulam theorem abound and their
applications include some of the more striking results in some fields (see
\cite{Mat94}). One of the more general formulations of Borsuk-Ulam type is
the theorem of Dold \cite{Dold}: For an $n$-connected $G$-space $X$ and $Y$
a free $G$-space of dimension at most $n$, there are \emph{no} $G$%
-equivariant mappings $X\rightarrow Y$. In this paper we consider a theorem
of this type for which the target space is the complement of an arrangement
of linear subspaces in a Euclidean space. Such spaces have been intensely
investigated in recent years and they represent natural test spaces for
problems in combinatorics and geometry. The nonexistence of an equivariant
mapping from a configuration space associated to a problem to a complement
of an arrangement means that equivariant mappings from the configuration
space to the Euclidean space containing the arrangement must meet the
arrangement, that is, the image must satisfy the linear conditions defining
the arrangment.

\noindent To control the algebraic topology of the complement of an
arrangement, we introduce the notion of a blow-up of a given arrangement
whose cohomology is especially nice \cite{Mark-Carsten}. The argument for
the main theorem is a novel use of the spectral sequence associated to the
Borel construction on a $G$-space.

\noindent The authors wish to acknowledge the hospitality of MSRI whose
atmosphere fosters collaboration. A great deal of gratitude goes to
professor C. Schultz for sharing his insight with us.

\subsection{Statement of the main result}

\noindent A finite family of linear subspaces $\mathcal{A}$ in some
Euclidean space $\mathbb{R}^{m}$ is known as an \textit{arrangement}. Let $%
M_{\mathcal{A}}$ denote the complement of the arrangment $\mathbb{R}%
^{m}\backslash \bigcup \mathcal{A}$. Suppose a group $G$ acts on $\mathbb{R}%
^{m}$. The set of fixed points of the action of $G$ is denoted by $(\mathbb{R%
}^{m})^{G}$. An arrangement $\mathcal{A}$ is a \textit{$G$-invariant
arrangement} if for all $g\in G$ and for all $L\in \mathcal{A}$, $gL\in
\mathcal{A}$. In the statement of the following theorem of Borsuk-Ulam type
we use notion of a blow up introduced in Section \ref{Sec:Blow-Up}.

\begin{theorem}
\label{Th:Main} Let $G$ denote a finite or a compact Lie group and $\Bbbk $
a field. Let $X$ be a $G$-space satisfying $H^{i}(X,\Bbbk )=0$ for $1\leq
i\leq n$ for some $n\geq 2$. Consider a $G$-invariant arrangement $\mathcal{A%
}$ in (some subspace $V $ of) $\mathbb{R}^{m}$ and its $G$-invariant blow up
$\mathfrak{B}(\mathcal{A})$ such that

\smallskip \noindent (A) the codimension of all maximal elements in $%
\mathcal{A}$ is less then $n+2$,

\noindent (B) for all $g\in G$ and all maximal elements $L$ of the
arrangement $\mathcal{A}$, we have $g\cdot L=L$,

\noindent (C) $G$ acts trivially on the cohomology $H^{\ast }(M_{\mathfrak{B}%
(\mathcal{A})},\Bbbk )$,

\noindent (D) the map $H^{\ast }(BG,\Bbbk )\rightarrow H^{\ast }(EG\times
_{G}M_{\mathfrak{B}(\mathcal{A})},\Bbbk )$, induced by the natural
projection $EG\times _{G}M_{\mathfrak{B}(\mathcal{A})}\rightarrow BG$, is
not a monomorphism, and

\noindent (E) for all $L\in \mathcal{A}$, $L\supseteq (\mathbb{R}^{m})^{G}$.

\smallskip\noindent Then there is no $G$-map $X\rightarrow M_{\mathcal{A}}$.
\end{theorem}

\begin{remark}
\textrm{Condition (A) implies that the complement $M_{\mathcal{A}}$ is $%
(n-1) $-connected. Theorem \ref{Th:Main} resembles Dold's theorem \cite[%
Remark on page 68]{Dold} in which condition (C) follows when the action on
the complement is free. The essential difference is not condition (C) (we
formulated it a little bit more generally) but that the dimension of the
space $M_{\mathcal{A}}$ can be arbitrary large. Also, it should be noted
that it is not generally possible to produce a $G$-invariant deformation of
the complement $M_{\mathcal{A}}$ to an $n$-dimensional subspace. Lemma 6.1
in \cite{BaMa2001} is only known example of an equivariant deformation of an
arrangement.}
\end{remark}

\begin{remark}
\textrm{Condition (A) can be substituted with the less restrictive condition
that
\begin{equation*}
\min \{\mathrm{codim}\ L\mid L\in \mathcal{A}\}=n+1.
\end{equation*}
Then in all the remaining conditions, replace the arrangement $\mathcal{A}$
by a subarrangement $\mathcal{A}^{\prime }$ generated by all maximal
elements of minimal codimension.}
\end{remark}

\begin{remark}
\textrm{Conditions (B) and (C) are in many cases equivalent. This can be
seen from the equivariant Goresky-MacPherson formula in \cite[Theorem
2.5.(ii)]{SunWelker}.}
\end{remark}

\begin{remark}
\textrm{Conditions (B) and (C) in some examples can be relaxed a little, but
not dropped all together. We illustrate this in Section \ref{sec:Conclusions}
with the construction of a $G$-map $X\rightarrow M_{\mathcal{A}}$ from an $n$%
-connected $G$-space $X$ to a codimension n+1 arrangement $\mathcal{A}$
complement. The arrangement $\mathcal{A}$ satisfies conditions (A), (D), (E)
but not (B) and (C). However, particular results can be obtained even when
the conditions (B) and (C) are not satisfied.}
\end{remark}

\begin{remark}
\textrm{When the group $G$ is connected, all arrangements satisfy condition
(C) since $\pi _{1}(BG)=\pi _{0}(G)=0$. Moreover, conditions (B) and (C) are
equivalent (from \cite[Theorem 2.5.(ii)]{SunWelker}). Condition (D) is
satisfied when, for example,}

\begin{itemize}
\item \textrm{$G$ acts freely on the complement $M_{\mathfrak{B}(\mathcal{A}%
)}$, or}

\item \textrm{$G$ is a $k$-torus or an elementary abelian $p$-group acting
without fixed points on $M_{\mathcal{A}}$ and consequently on $M_{\mathfrak{B%
}(\mathcal{A})}$.}
\end{itemize}
\end{remark}

\subsection{\label{Sec:appResults}An Application: Antipodal chees problem}

\noindent Every theorem of general type is usually a product of successful
or more often unsuccessful effort in solving some concrete problem. One of
the motivating problems for the study of Borsuk Ulam type theorems for
complements of arrangements is a class of mass partition problems discussed
in \cite{BaMa2001}, \cite{BaMa2002}, \cite{Bl}, \cite{Bl-Di} and \cite%
{Bl-Vr-Ziv}. Particularly, we present a problem which solution, after
careful combinatorial reformulation, presents a direct consequence of
Theorem \ref{Th:Main}.

\subsubsection*{Antipodal cheese problem}

\begin{figure}[tbh]
\centering
\includegraphics[scale=0.80]{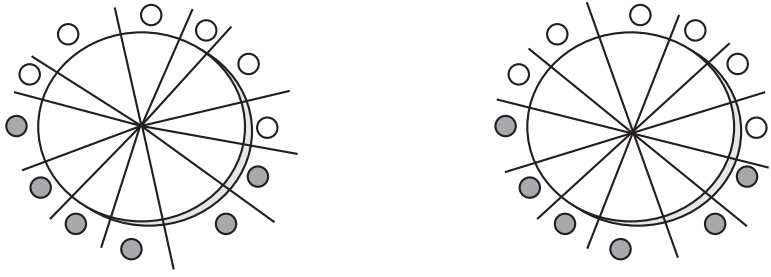}
\caption{{\protect\small Half-straight and straight cheese cuts.}}
\label{Fig:Cheese}
\end{figure}

\noindent Suppose that even number $2k$ of people is sitting around a circle
table in such a way that every one has its antipodal friend. On the table
there is pile of $j$ (high dimensional) cheese pieces (in $%
\mathbb{R}
^{d}$); all of different shape, mass, density and flavor. The line knife is
available and the cheese can be cut only simulataneously, all $j$ peaces at
once. There are two types of cuts we allow (Figure \ref{Fig:Cheese}):

\begin{itemize}
\item half-straight cut: pick a point on the table and make $2k$ straight
(line) cuts starting at the chosen point and continuing in one direction,

\item straight cut: pick a point on the table and make $k$ straight (line)
cuts through the chosen point in both directions.
\end{itemize}

\noindent The objective is to divide cheese ($j$ pieces) in $%
\mathbb{R}
^{d}$ by $2k$ half-straight cuts or $k$ straight cuts in such a way that
every member of an antipodal pear get the same, non-negative, part of each
piece of cheese.

\subsubsection*{Mathematical reformulation}

\noindent The vocabulary for mathematical translation of the antipodal
cheese problem is:
\begin{eqnarray*}
\text{half-straight cut } &\rightarrow &\text{ fan,} \\
\text{straight cut } &\rightarrow &\text{ arrangement in fan position,} \\
\text{piece of cheese } &\rightarrow &\text{ measure.}
\end{eqnarray*}

\medskip

\noindent Let $H$ be a hyperplane in $\mathbb{R}^{d}$ and $L$ a codimension
one subspace inside $H$. The connected components of a space $H\backslash L$
are \textbf{half-hyperplanes} determined by a pair $(H,L)$. If $F_{1}$ and $%
F_{2}$ are the half-hyperplane determined by the pair $(H,L)$, then $L$ is
the boundary of both half-hyperplanes.

\begin{definition}
A $k$\textbf{-fan} in $\mathbb{R}^{d}$ is a collection $(L;F_{1},\ldots,
F_{k})$ consisting of

\noindent (A) a $(d-2)$-dimensional oriented linear subspace $L$, and

\noindent (B) different half-hyperplanes $F_{1},\ldots, F_{k}$ with the
common boundary $L$, oriented by a compatible orientation on the plane $%
L^{\bot }$.
\end{definition}

\noindent Let $S(L^{\perp })$ denote the unit circle lying in the plane $%
L^{\perp }$. Condition 7.(B) suggests that the intersection points $%
F_{1}\cap S(L^{\bot }),\ldots ,F_{k}\cap S(L^{\bot })$ are consecutive
points on the circle $S(L^{\bot })$ oriented consistently with the given
orientation on $L^{\bot }$. The $k$\textbf{-fan }$(L,l_{1},\ldots ,l_{k})$
on the sphere $S^{d-1}\subset \mathbb{R}^{d}$ is the trace of a $k$-fan $%
(L;F_{1},\ldots ,F_{k})$ in $\mathbb{R}^{d}$ obtained by slicing the sphere $%
S^{d-1}$ along half-hyperplanes, $l_{i}=S^{d-1}\cap F_{i}$. Sometimes,
instead of a sequence $l_{1},\ldots ,l_{k}$ of cuts for the model of a fan,
we will prefer the sequence of open sets, called \emph{orthants}, $\mathcal{O%
}_{i}$, $i\in \{1,\ldots ,k\}$ on the sphere $S^{d-1}$ lying between
consecutive cuts $l_{i},l_{i+1}$, $i\in \{1,\ldots ,k\}$ (we assume $%
l_{k+1}=l_{1}$). A third model for a $k$-fan is the collection $%
(L;v_{1},\ldots ,v_{k})$ of a codimension 2 subspace $L$ inside $\mathbb{R}%
^{d}$ and $k$ vectors $v_{1},\ldots ,v_{k}$ on the circle $S(L^{\bot })\cong
S^{1}$. Denote by $\phi _{i}$ angle between $v_{i}$ and $v_{i+1}$ ($%
v_{k+1}=v_{1}$). Then $\phi _{1}+\cdots +\phi _{k}=2\pi $. The space of all $%
k$-fans in $\mathbb{R}^{d}$ or on the sphere $S^{d-1}$ is denoted by $%
\mathcal{F}_{k}$.

\begin{definition}
A coordinate hyperplane arrangement $\mathcal{A}=\{H_{1},\ldots, H_{k}\}$ is
\textbf{in $k$-fan position} if the intersection $H_{1}\cap \cdots \cap
H_{k} $ is a subspace of codimension one inside each $H_{i}$. In other
words, $\mathcal{A}$ is in a fan position if there is a $2k$-fan $%
(L;F_{1},\ldots, F_{2k}) $ such that $H_{1}\cap \cdots \cap H_{k}=L$ and $%
H_{1}\cup \cdots \cup H_{k}=L\cup F_{1}\cup \cdots \cup F_{2k}$.
\end{definition}

\begin{figure}[tbh]
\centering\includegraphics[scale=0.80]{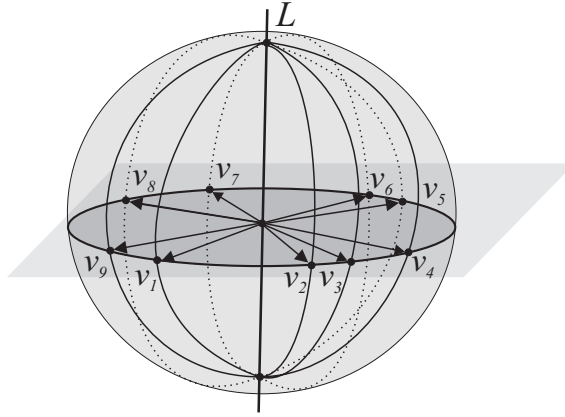}
\caption{{A $9$-fan in $\mathbb{R}^{3}$}}
\end{figure}

\noindent As in the case of a fan, an arrangement $\mathcal{A}%
=\{H_{1},\ldots ,H_{k}\}$ in fan position inherits a natural orientation
from $L^{\bot }$, $L=H_{1}\cap \cdots \cap H_{k}$. Besides the natural
ordering of hyperplanes, the orientation on $L^{\bot }$ induces an
orientation of the connected components (orthants) of the complement $%
\mathbb{R}^{d}\backslash \bigcup H_{i}$. The orientation is determined up to
a cyclic permutation. If $(H_{1},\ldots ,H_{k})$ is the induced ordering and
$H_{k+1}=H_{1}$, then we denote by

\begin{itemize}
\item $\mathcal{O}_{i}^{+}$ the orthant between $H_{i}$ and $H_{i+1}$, and by

\item $\mathcal{O}_{i}^{-}$ the orthant between $H_{i+1}$ and $H_{i}$.
\end{itemize}

\bigskip

\noindent A measure $\mu $ on a sphere $S^{d-1}$ is a \textbf{proper}
measure if for every hyperplane $H\subset \mathbb{R}^{d}$, $\mu (H\cap
S^{d-1})=0$ and for every non-empty open set $U\subseteq S^{d-1}$, $\mu
(U)>0 $. From now on a measure on a sphere $S^{d-1}$ will mean a proper
Borel probability measure.

\bigskip

\noindent Rational vector $\mathbf{\alpha }=(\alpha _{1},\ldots ,\alpha
_{k},\alpha _{k+1},\ldots ,\alpha _{2k})\in \mathbb{Q}^{2k}$ satisfying
\begin{equation*}
\text{for all $i\in \{1,\ldots ,k\}$},\quad \alpha _{i}>0,\,\alpha
_{i}=\alpha _{k+i},
\end{equation*}%
and%
\begin{equation*}
\sum_{i=1}^{k}\alpha _{i}=\tfrac{1}{2}
\end{equation*}%
is called\emph{\ }\textbf{ration}.

\begin{definition}
Let $\mathcal{M}=\{\mu _{1},\ldots ,\mu _{j}\}$ be a collection of measures
on $S^{d-1}$.
A $2k$-fan $(L,\mathcal{O}_{1},\ldots ,\mathcal{O}_{2k})$ is a $\mathbf{%
\alpha }$\textbf{-partition} of $\mathcal{M}$ if there is $i\in \{1,\ldots
,j\}$ such that%
\begin{equation*}
\text{for all $t\in \{1,\ldots ,2k\}$, }\mu _{i}(\mathcal{O}_{t})=\alpha _{t}
\end{equation*}%
and%
\begin{equation*}
\text{for all $t\in \{1,\ldots ,k\}$ and all $l\in \{1,\ldots ,j\}$},\ \mu
_{l}(\mathcal{O}_{t})=\mu _{l}(\mathcal{O}_{k+t})\text{.}
\end{equation*}

\noindent An arrangement $\mathcal{A}=\{H_{1},\ldots ,H_{k}\}$ in $k$-fan
position is an $\mathbf{\alpha }$\textbf{-partition} of $\mathcal{M}$ if
there is $i\in \{1,\ldots ,j\}$ such that%
\begin{equation*}
\text{for all $t\in \{1,\ldots ,2k\}$, }\mu _{i}(\mathcal{O}_{t})=\alpha _{t}
\end{equation*}%
and%
\begin{equation*}
\text{for all $t\in \{1,\ldots ,k\}$ and all $l\in \{1,\ldots ,j\}$},\ \mu
_{l}(\mathcal{O}_{t})=\mu _{l}(\mathcal{O}_{k+t})\text{.}
\end{equation*}
\end{definition}

\subsubsection*{A consequence of Theorem \protect\ref{Th:Main}}

\noindent The solutions of anitpodal cheese problem implied by Theorem \ref%
{Th:Main} are given in the following theorem.

\begin{theorem}
\label{Th:APP}\qquad \newline
(A) Let $k>0$ be an integer and $\mathbf{\alpha }\in \mathbb{Q}^{2k}$ be a
ration. If $k(j-1)<d-1$, then for every collection of $j$ measures $\mathcal{%
M}$ on a sphere $S^{d-1}$, there exists a $\mathbf{\alpha }$-partition of $%
\mathcal{M}$ by a $k$-fan

\noindent (B) Let $k>0$ be an integer and $\mathbf{\alpha }\in \mathbb{Q}%
^{2k}$ be a ration. If $kj<d-1$, then for every collection of $j$ measures $%
\mathcal{M}$ on a sphere $S^{d-1}$, there exists an $\mathbf{\alpha }$%
-partition of $\mathcal{M}$ by an arrangement in a $k$-fan position.
\end{theorem}

\section{Proof of Theorem \protect\ref{Th:Main}}

\noindent The proof is carried out using the Borel construction and its
associated Serre spectral sequence. Given an equivariant mapping between $G$%
-spaces, $f\colon X\rightarrow Y$, there is an induced mapping of Borel
constructions, $EG\times _{G}f\colon EG\times _{G}X\rightarrow EG\times
_{G}Y $. The notion of a blow up of an arrangement is the key construction
which leads to the proof of the main theorem.

\subsection{\label{Sec:Blow-Up}Blow up of an arrangement}

\noindent By the \textbf{codimension} of an arrangement $\mathcal{A}$,
denoted \textrm{codim}$_{\mathbb{R}^{m}}\mathcal{A}$, we understand
\begin{equation*}
\mathrm{codim}_{\mathbb{R}^{m}}\mathcal{A}=\min_{L\in \mathcal{A}}\left\{
\mathrm{codim}_{\mathbb{R}^{m}}L\right\} .
\end{equation*}%
Following Definition 5.3 in \cite{Mark-Carsten}, an arrangement $\mathcal{A}$
is a $c$\emph{-\textbf{arrangement} }if

\begin{itemize}
\item for every maximal element $L$ in $\mathcal{A}$, $\mathrm{codim}_{%
\mathbb{R}^{m}}L=c$

\item for all pairs $L_{1}\subset L_{2}$ of elements in $\mathcal{A}$, $c$
divides $\mathrm{codim}_{L_{2}}L_{1}$.
\end{itemize}

\medskip

\noindent Recall that if $X$ and $Y$ are $G$-spaces, the \emph{diagonal
action} of $G$ on the product $X\times Y$ is given by $g\cdot (x,y)=(g\cdot
x,g\cdot y)$. By $G\mathcal{A}$ we denote the \emph{minimal }$G$\emph{%
-invariant} arrangement containing the arrangement $\mathcal{A}$, namely, $G%
\mathcal{A}=\{gL\mid g\in G\text{ and }L\in \mathcal{A}\}$. An arrangement $%
\mathcal{A}$ is $G$-invariant if and only if $G\mathcal{A=A}$ .

If $L\subset \mathbb{R}^{m}$ is a linear subspace, recall that there exists
a family of \emph{forms}, $\xi _{1}$, \dots , $\xi _{t}$, given by
\begin{equation*}
\xi _{j}(x_{1},\ldots ,x_{m})=a_{1j}x_{1}+\cdots +a_{mj}x_{m}
\end{equation*}%
for which $L=\{\mathbf{x}\in \mathbb{R}^{m}\mid \xi _{1}(\mathbf{x})=\cdots
=\xi _{s}(\mathbf{x})=0\}$.

\begin{definition}
Let $\mathcal{A}$ be an arrangement of linear subspaces in $\mathbb{R}^{m}$,
$\{L_{1},\ldots, L_{w}\}$ the set of maximal elements of $\mathcal{A}$ and $%
k_{i}=\mathrm{codim}_{\mathbb{R}^{m}}L_{i}$ , for $i\in \{1,\ldots, w\}$.
For each maximal element $L_{i}$, there is a family $\{\xi _{i,1},\ldots,
\xi _{i,k_{i}}\}$ of (linearly independent) forms defining $L_{i}$. The
\textbf{blow up} of the arrangement $\mathcal{A}$ is the arrangement $%
\mathfrak{B}(\mathcal{A})$ in
\begin{equation*}
(\mathbb{R}^{m})^{k_{1}+\cdots +k_{w}}=( \mathbb{R}^{m})^{k_{1}}\times
\cdots \times ( \mathbb{R}^{m})^{k_{w}}=E_{1}\times \cdots \times E_{w}
\end{equation*}%
(where $\left( \mathbb{R}^{m}\right)^{k_{i}}=E_{i}$) is defined by $w$
maximal elements $\bar{L}_{1}, \ldots,\bar{L}_{w}$ introduced in the
following way. The subspace $\bar{L}_{i}$, $i=1,\ldots, w$, is defined by
forms:

$\xi _{i,1} = 0$ seen as a form on the $1$-st copy of $\mathbb{R}^{n}$ in $%
E_{i}$;

$\xi _{i,2} = 0$ seen as a form on the $2$-nd copy of $\mathbb{R}^{n}$ in $%
E_{i}$;

$\vdots$

$\xi _{i,k_{i}} = 0$ seen as a form on the $k_{i}$-th copy of $\mathbb{R}%
^{n} $ in $E_{i}$.
\end{definition}

\noindent The blow up $\mathfrak{B}(\mathcal{A})$ depends on the choice of
the linear forms $\xi _{\ast ,\ast }$. Observe that we do not allow any
extra dependent forms. Note also that the arrangement operations $\mathfrak{%
B(\cdot )}$ and $G(\cdot )$ do not commute.

\begin{remark}
\textrm{For an arrangement $\mathcal{A}$ inside (an invariant $G$-) subspace
$V\subset \mathbb{R}^{m}$, the blow up is an arrangement inside $(V)
^{k_{1}}\times \cdots \times ( V) ^{k_{w}}$ defined as in Definition~11.}
\end{remark}

\begin{example}
\textrm{Let $L\subset \mathbb{R}^{2}$ denote the trivial subspace $L =
\{(0,0)\}$, and $\mathcal{A}=\{L\}$. Then the blow up $\mathfrak{B}(\mathcal{%
A})$ is an arrangement in $\mathbb{R}^{4}$ with one element defined by $%
x_{1}=x_{4}=0$.}
\end{example}

\noindent Here is a list of significant properties of the blow up of
arrangement.

\begin{proposition}
\label{prop:blow-Up-1}Let $\mathcal{A}$ be an arrangement of linear
subspaces in $\mathbb{R}^{m}$ and $\mathfrak{B}(\mathcal{A})$ its associated
blow up in $( \mathbb{R}^{m}) ^{k_{1}+\cdots +k_{w}}$.

\noindent (A) $\mathrm{codim}_{\mathbb{R}^{m}}\mathcal{A=}\mathrm{codim}_{%
\mathbb{R}^{m(k_{1}+\cdots +k_{w})}}\mathfrak{B}(\mathcal{A})$.

\noindent (B) $\mathfrak{B}(\mathcal{A})$ is a $\left( \mathrm{codim}_{%
\mathbb{R}^{m}}\mathcal{A}\right) $-arrangement.

\noindent (C) The identity map $\mathbb{R}^{m}\rightarrow \mathbb{R}^{m}$
induces the diagonal map $D\colon \mathbb{R}^{m}\rightarrow \left( \mathbb{R}%
^{m}\right) ^{k_{1}+\cdots +k_{w}}$ which restricts to a map of complements%
\begin{equation*}
D\colon \mathbb{R}^{m}\backslash \bigcup \mathcal{A~\rightarrow ~}\left(
\mathbb{R}^{m}\right) ^{k_{1}+\cdots +k_{w}}\backslash \bigcup \mathfrak{B}(%
\mathcal{A}).
\end{equation*}
\end{proposition}

\begin{proof}
These statements are direct consequences of the blow up construction.
\end{proof}

\begin{proposition}
\label{prop:blow-Up-2} Consider a $G$-action on $\mathbb{R}^{m}$, which
extends diagonally to the product $\left( \mathbb{R}^{m}\right)
^{k_{1}+\cdots +k_{w}}$. Let $\mathcal{A}$ be a $G$-invariant arrangement in
$\mathbb{R}^{m}$ such that conditions (B) and (C) of Theorem \ref{Th:Main}
are satisfied. Then the defining forms can be chosen in such a way that $%
\mathfrak{B}(\mathcal{A})$ is also a $G$-invariant arrangement satisfying
the same conditions, and more:

\noindent (A) The diagonal map $D\colon \mathbb{R}^{m}\rightarrow (\mathbb{R}%
^{m})^{k_{1}+\cdots +k_{w}}$ is a $G$-map. Moreover, the diagonal map
restricts to a $G$-map of complements%
\begin{equation*}
D\colon \mathbb{R}^{m}\backslash \bigcup \mathcal{A}=M_{\mathcal{A}}\mathcal{%
~\rightarrow ~}(\mathbb{R}^{m})^{k_{1}+\cdots +k_{w}}\backslash \bigcup
\mathfrak{B}(\mathcal{A})=M_{\mathfrak{B}(\mathcal{A})};
\end{equation*}

\noindent (B) If $k_{1}=\cdots =k_{w}=k$, then the blow up $\mathfrak{B}(%
\mathcal{A})$ is a $k$-arrangement and the cohomology ring $\tilde{H}^{\ast
}(M_{\mathfrak{B}(\mathcal{A})},\Bbbk )$ is generated as an algebra by ${H}%
^{k-1}(M_{\mathfrak{B}(\mathcal{A})},\Bbbk )$.

\noindent (C) If for all $L\in \mathcal{A}$, $L\supseteq (\mathbb{R}%
^{m})^{G} $, then, for all $L\in \mathfrak{B}(\mathcal{A})$, we have $%
L\supseteq \left( (\mathbb{R}^{m})^{k_{1}+\cdots +k_{w}}\right) ^{G}.$

\noindent (D) If the map $H^{\ast }(BG,\Bbbk )\rightarrow H^{\ast }(EG\times
_{G}M_{\mathfrak{B}(\mathcal{A})},\Bbbk )$ is not a monomorphism, then the
same will be true for the map $H^{\ast }(BG,\Bbbk )\rightarrow H^{\ast
}(EG\times _{G}M_{\mathcal{A}},\Bbbk )$.
\end{proposition}

\begin{proof}
Let $\mathcal{A}$ be a $G$-invariant arrangement satisfying conditions (B)
and (C) of Theorem \ref{Th:Main}. For $L$ a maximal element of $\mathcal{A}$%
, we can choose defining forms $\{\xi _{i,1},\ldots, \xi _{i,k}\}$ in such a
way that for all $g\in G$ and all $j\in \{1,\ldots, k\}$ holds $g\cdot \xi
_{i,j}=\pm \xi _{i,j}$. With this choice of defining forms and because we
assumed condition (C) of Theorem \ref{Th:Main} the blow up $\mathfrak{B}(%
\mathcal{A}) $ satisfies the required properties. Statement (A) is a
consequence of the diagonal action on $(\mathbb{R}^{m})^{k_{1}+\cdots
+k_{w}} $. The first part of the property (B) follows from assuming
condition (B) of Theorem \ref{Th:Main}. The second part of (B) follows from
Corollary 5.6 in \cite{Mark-Carsten}. The equality $\left( (\mathbb{R}%
^{m})^{k_{1}+\cdots +k_{w}}\right) ^{G}=((\mathbb{R}^{m})^{G})^{k_{1}+\cdots
+k_{w}}$ implies (C).

To prove Statement (D) we consider the mapping induced by the $G$%
-equivariant diagonal mapping $M_{\mathcal{A}}\rightarrow M_{\mathfrak{B}(%
\mathcal{A})}$ on the Borel constructions,
\begin{equation*}
D\colon EG\times _{G}M_{\mathcal{A}}\rightarrow EG\times _{G}M_{\mathfrak{B}(%
\mathcal{A})}.
\end{equation*}%
By assumption, the edge homomorphism $H^{\ast }(BG,\Bbbk )\rightarrow
H^{\ast }(EG\times _{G}M_{\mathfrak{B}(\mathcal{A})},\Bbbk )$ is not a
monomorphism and this is equivalent to the fact that there is a nonzero
differential in the Serre spectral sequence for $EG\times _{G}M_{\mathfrak{B}%
(\mathcal{A})}\rightarrow EG\times_{G} \{\hbox{\rm pt}\} =BG$. By
assumption, the $E_{2}$-term may be written $E_{2}^{p,q}\cong H^{p}(BG,\Bbbk
)\otimes H^{q}(M_{\mathfrak{B}(\mathcal{A})},\Bbbk )$. By property (B), $%
\tilde{H}^{\ast }(M_{\mathfrak{B}(\mathcal{A})},\Bbbk )$ is generated as an
algebra in dimension $k-1$. Since the cohomology Serre spectral sequence is
multiplicative, the first differential must be $d_{k}\colon \tilde{H}%
^{k-1}(M_{\mathfrak{B}(\mathcal{A})},\Bbbk )\rightarrow H^{k}(BG,\Bbbk )$.
If $d_{k}=0$, then $d_{k+l}=0$ for $l\geq 1$, and the spectral sequence
collapses at $E_{2}$, which contradicts the assumption that the edge
homomorphism is not a monomorphism.

Suppose $1\otimes v=d_{k}(u\otimes 1)$ for $v\in H^{\ast }(BG,\Bbbk )$ and $%
u\in H^{\ast }(M_{\mathfrak{B}(\mathcal{A})},\Bbbk )$. The diagonal mapping
induces a mapping of spectral sequences that is given on the $E_{2}$-term by
the identity on $E_{2}^{\ast ,0}$ and the induced mapping on cohomology on $%
E_{2}^{0,\ast }$. Since the differential commutes with this induced mapping,
we have
\begin{equation*}
0\neq 1\otimes v=E_{2}(D)(1\otimes v)=E_{2}(D)(d_{k}(u\otimes
1))=d_{k}(D^{\ast }(u)\otimes 1).
\end{equation*}%
If $D^{\ast }(u)=0$, then $d_{k}(D^{\ast }(u)\otimes 1)=0$, which
contradicts the choice of $v$. Thus, $D^{\ast }(u)\neq 0$ and the
differential $d_{k}\neq 0$ on the spectral sequence for $EG\times _{G}M_{%
\mathcal{A}}\rightarrow BG$. Statement (C) follows immediately.
\end{proof}

\subsection{Comparing Serre spectral sequences; proof of Theorem \protect\ref%
{Th:Main}}

\noindent For simplicity reasons let us assume that the codimension of all
maximal elements of the arrangement $\mathcal{A}$ is $n+1$.

\medskip

\noindent To prove Theorem \ref{Th:Main}\ we made assumption (D), that the
mapping
\begin{equation*}
H^{\ast }(BG,\Bbbk )\rightarrow H^{\ast }(EG\times _{G}M_{\mathfrak{B}(%
\mathcal{A})},\Bbbk )
\end{equation*}%
is not a monomorphism. By the assumption (C) and the choice of a field $%
\Bbbk $ for coefficients, we can write the $E_{2}$-term of the spectral
sequence for the Borel construction $EG\times _{G}M_{\mathfrak{B}(\mathcal{A}%
)}$ as
\begin{equation*}
E_{2}^{\ast ,\ast }\cong H^{\ast }(BG,\Bbbk )\otimes H^{\ast }(M_{\mathfrak{B%
}(\mathcal{A})},\Bbbk ).
\end{equation*}%
It follows from the assumption (A) that $M_{\mathfrak{B}(\mathcal{A})}$ is $%
n $-connected, and so there is a nonzero differential $d_{n+1}\colon
E_{n+1}^{0,n}\rightarrow E_{n+1}^{n+1,0}$ in this spectral sequence, as
argued in the proof of Proposition \ref{prop:blow-Up-2}. We apply this
observation to study the existence of a $G$-map $f\colon X%
\rightarrow M_{\mathcal{A}}$. As in Proposition \ref{prop:blow-Up-2}, we
have the $G$-map $D\colon M_{\mathcal{A}}\rightarrow M_{\mathfrak{B}(%
\mathcal{A})}$ and consequently a $G$-map%
\begin{equation}
X\overset{f}{\rightarrow }M_{\mathcal{A}}\overset{D}{\rightarrow }M_{%
\mathfrak{B}(\mathcal{A})}.  \label{map1}
\end{equation}

\noindent Then there is an induced homomorphism $(f\circ D)^{\ast }$ in
equivariant cohomology:%
\begin{equation*}
H^{\ast }(X\times _{G}EG;\Bbbk )\overset{f^{\ast }}{\leftarrow }H^{\ast }(M_{%
\mathcal{A}}\times _{G}EG;\Bbbk )\overset{D^{\ast }}{\leftarrow }H^{\ast
}(M_{\mathfrak{B}(\mathcal{A})}\times _{G}EG;\Bbbk )\text{,}
\end{equation*}%
as well as an induced map of Serre spectral sequences
\begin{equation*}
E_{\ast }^{\ast .\ast }(X\times _{G}EG;\Bbbk )\overset{E^{\ast }(f\circ D)}{%
\leftarrow }E_{\ast }^{\ast ,\ast }(M_{\mathfrak{B}(\mathcal{A})}\times
_{G}EG;\Bbbk )
\end{equation*}%
associated to fibrations%
\begin{equation*}
\begin{array}{ccccccccccc}
X & \rightarrow & X\times _{G}EG &  &  &  &  &  & M_{\mathfrak{B}(\mathcal{A}%
)} & \rightarrow & M_{\mathfrak{B}(\mathcal{A})}\times _{G}EG \\
&  & \downarrow &  &  &  &  &  &  &  & \downarrow \\
&  & BG &  &  &  &  &  &  &  & BG%
\end{array}%
\end{equation*}

\noindent We use the naturality of the spectral sequence and the property
that $E_{2}^{\ast ,0}(f\circ D)=\mathrm{id}$ on $H^{\ast }(BG,\Bbbk )$. The
induced map $E_{2}(f\circ D)$ and the $E_{2}$-term of the associated
spectral sequences for $EG\times _{G}X$ and $EG\times _{G}M_{\mathfrak{B}(%
\mathcal{A})}$ can be pictured:

\begin{figure}[tbh]
\centering
\includegraphics[scale=0.60]{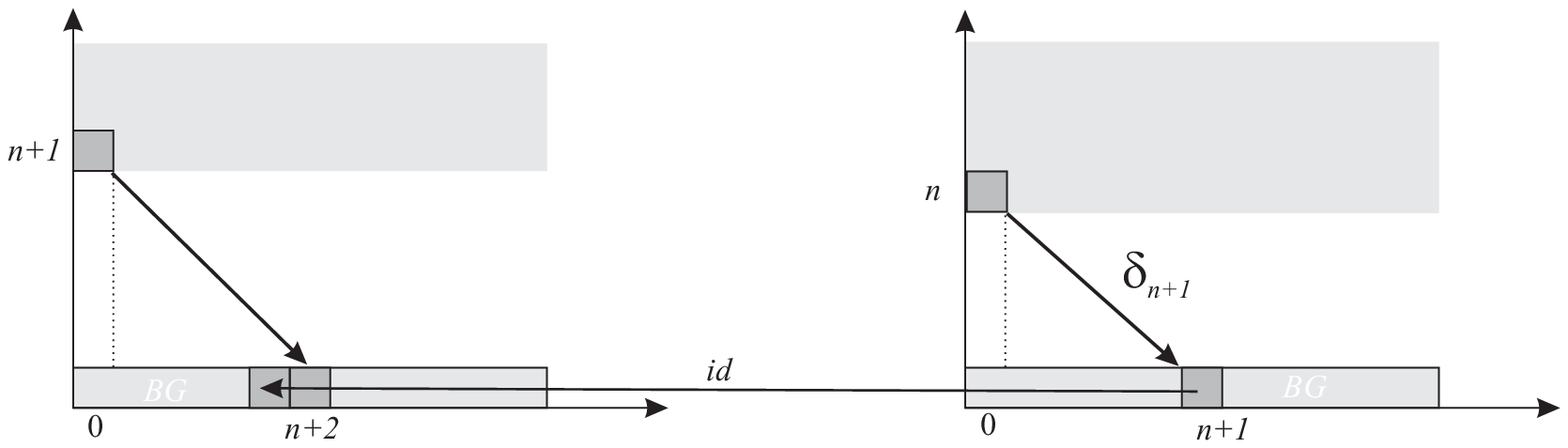}
\caption{}
\label{Fig:tableu2}
\end{figure}

\noindent The first nonzero differential is $d_{n+1}\colon
E_{n+1}^{0,n}=E_{2}^{0,n}\rightarrow E_{n+1}^{n+1,0}=E_{2}^{n+1,0}$ in the
spectral sequence for $EG\times _{G}M_{\mathfrak{B}(\mathcal{A})}$. By the
same argument in Proposition \ref{prop:blow-Up-2}, there must be a nonzero
class in $H^{n}(X,\Bbbk )$ to be the image of a generator $u\in H^{n}(M_{%
\mathfrak{B}(\mathcal{A})},\Bbbk )$ which transgresses to $H^{n+1}(BG,\Bbbk
) $. However, the connectivity of $X$ presents no nonzero classes in this
dimension, nor any below this dimension from which to launch a differential.
Thus, the existence of the $G$-map leads to a contradiction.

\noindent The proof of the main theorem is now complete.

\section{\label{sec:Conclusions}Is it possible to obtain more?}

\noindent The following examples will illustrate that Theorem \ref{Th:Main}
is the best general theorem of Borsuk-Ulam type one can obtain for
complements of arrangements one can obtain.

\subsection{\label{Sec:Example}An example}

\noindent Let $G$ be the cyclic group $\mathbb{Z}/{n}=\langle \omega \rangle
$ where $n>2$ is odd. Let $G$ act \textbf{freely} on $S^{3}$ and on $\mathbb{%
R}^{n}$ by the cyclic shift, $\omega \cdot (x_{1},\ldots
,x_{n})=(x_{2},\ldots ,x_{n},x_{1})$. Let $W_{n}$ denote the $G$-invariant
subspace $\{(x_{1},\ldots ,x_{n})~|~x_{1}+\cdots +x_{n}=0\}$. Consider the
minimal $G$-invariant arrangement $\mathcal{A}$ containing the subspace $%
L\subset W_{n}^{\oplus 3}$ defined by%
\begin{equation*}
L=\{({\mathbf{x}}_{i,1};{\mathbf{x}}_{i,2};{\mathbf{x}}_{i,3})_{i\in
\{1,\ldots ,n\}}~|~{\mathbf{x}}_{1,1}={\mathbf{x}}_{1,2}={\mathbf{x}}%
_{1,3}=0\}.
\end{equation*}%
Endow $W_{n}^{\oplus 3}$ with the diagonal $G$-action. \emph{Then there is a
}$G$\emph{-map }$S^{3}\rightarrow W_{n}^{\oplus 3}\backslash \bigcup A=M_{%
\mathcal{A}}$\emph{\ and the }$G$-\emph{action on }$H^{\ast }(M_{\mathcal{A}%
},\Bbbk )$\emph{\ is nontrivial for any field }$\Bbbk $\emph{.}

\medskip

\noindent To prove this statement we use equivariant obstruction theory (for
background, consult \cite{Dieck87} and for applications \cite{Bl-Di}). The
obstruction theory can be applied because the action of $\mathbb{Z}/{n}$ on
the sphere $S^{3}$ is assumed to be free. Since $\mathrm{codim}%
_{W_{n}^{\oplus 3}}\mathcal{A=}3$ the complement $M_{\mathcal{A}}$ is $1$%
-connected, $2$-simple and therefore the problem of the existence of a map $%
S^{3}\rightarrow M_{\mathcal{A}}$ depends on the primary obstruction. The
primary obstruction lives in $H_{G}^{3}(S^{3},\pi _{2}(M_{\mathcal{A}}%
\mathcal{))}$.

\begin{lemma}
$H_{G}^{3}(S^{3},\pi _{2}(M_{\mathcal{A}}\mathcal{))\cong\ }\mathbb{Z}$.
\end{lemma}

\begin{proof}
The complement $M_{\mathcal{A}}$ is $1$-connected and according to the
Hurewicz theorem
\begin{equation*}
\pi _{2}(M_{\mathcal{A}})=H_{2}(M_{\mathcal{A}},\mathbb{Z}).
\end{equation*}%
Since $\mathcal{A}$ is a 3-arrangement, the Goresky-MacPherson formula
applied to $\mathcal{A}$ implies that $H_{2}(M_{\mathcal{A}},\mathbb{Z})=%
\mathbb{Z}[\mathbb{Z}/{n}]$ as a $\mathbb{Z}/{n}$-module. The equivariant
Poincar\'{e} duality isomorphism \cite[Theorem 1.4]{Bl-Di} applied to the $G$%
-manifold $S^{3}$ yields an isomorphism
\begin{equation*}
H_{G}^{3}(S^{3},H_{2}(M_{\mathcal{A}},\mathbb{Z}))\cong
H_{0}^{G}(S^{3},H_{2}(M_{\mathcal{A}},\mathbb{Z})\otimes \mathcal{Z})
\end{equation*}%
where $\mathcal{Z}$ is the $G$-module $H_{n+1}(S^{3},\mathbb{Z})\cong
\mathbb{Z}$. Since $n$ is odd, the $G$-module $\mathcal{Z}$ is trivial and
therefore,%
\begin{equation*}
H_{G}^{3}(S^{3},H_{2}(M_{\mathcal{A}},\mathbb{Z}))\cong
H_{0}^{G}(S^{3},H_{2}(M_{\mathcal{A}},\mathbb{Z})).
\end{equation*}%
From homological algebra (for example, \cite[(1.5), p.57]{Brown}) we have
that%
\begin{equation*}
H_{0}^{G}(S^{3},H_{2}(M_{\mathcal{A}},\mathbb{Z}))\cong H_{0}^{G}(G,H_{2}(M_{%
\mathcal{A}},\mathbb{Z}))\cong H_{2}(M_{\mathcal{A}},\mathbb{Z})_{G}.
\end{equation*}%
Thus $H_{G}^{3}(S^{3},H_{2}(M_{\mathcal{A}},\mathbb{Z}))\cong \left( \mathbb{%
Z}[\mathbb{Z}/{n}]\right) _{\mathbb{Z}/{n}}\cong \mathbb{Z}$. For general
discussion of this type of argument see \cite[Section 1]{Bl-Di}.
\end{proof}

\begin{lemma}
The primary obstruction is a torsion element in $H_{G}^{3}(S^{3},\pi _{2}(M_{%
\mathcal{A}}\mathcal{))}$.
\end{lemma}

\begin{proof}
Let $H$ be a subgroup of $G$. There is a natural restriction\ map
\begin{equation*}
r:H_{G}^{3}(S^{3},H_{2}(M_{\mathcal{A}},\mathcal{\mathbb{Z})})\rightarrow
H_{H}^{3}(S^{3},H_{2}(M_{\mathcal{A}},\mathcal{\mathbb{Z})}),
\end{equation*}%
which on the cochain level is just the forgetful map sending a $G$-cochain
to the same cochain interpreted as a $H$-cochain. From the geometric
definition of the obstruction cocycle \cite[Section 1.2]{Bl-Di} it follows
that the restriction of the primary obstruction is the primary obstruction
for the existence of a $H$-map $H_{2}(M_{\mathcal{A}},\mathcal{\mathbb{Z})}$%
. Furthermore there exists a natural map%
\begin{equation*}
\tau :H_{H}^{3}(S^{3},H_{2}(M_{\mathcal{A}},\mathcal{\mathbb{Z})}%
)\rightarrow H_{G}^{3}(S^{3},H_{2}(M_{\mathcal{A}},\mathcal{\mathbb{Z})}).
\end{equation*}%
in the opposite direction given by the transfer. It is known \cite[Section
III.9. Proposition 9.5.(ii)]{Brown} that the composition of the restriction
with the transfer is just multiplication by the index $[G:H]$:
\begin{equation*}
\begin{array}{ccccc}
H_{G}^{3}(S^{3},H_{2}(M_{\mathcal{A}},\mathcal{\mathbb{Z})}) & \overset{r}{%
\rightarrow } & H_{H}^{3}(S^{3},H_{2}(M_{\mathcal{A}},\mathcal{\mathbb{Z})})
& \overset{\tau }{\rightarrow } & H_{G}^{3}(S^{3},H_{2}(M_{\mathcal{A}},%
\mathcal{\mathbb{Z})}).%
\end{array}%
\end{equation*}%
Consider $H=\{e\}$, the trivial group. The $H$ primary obstruction is zero,
since there is an $H$-map $S^{3}\rightarrow M_{\mathcal{A}}$. Therefore the
primary $G$-obstruction multiplied by the index $[G:H]$ vanishes, that is,
the primary obstruction is a torsion element in $H_{G}^{3}(S^{3},H_{2}(M_{%
\mathcal{A}},\mathcal{\mathbb{Z})})$.
\end{proof}

\smallskip \noindent The primary obstruction responsible for the existence
of a $G$-map, $S^{3}\rightarrow M_{\mathcal{A}}$, as a torsion element in $%
\mathbb{Z} $, must vanish. Therefore, there exists a $G$-map $%
S^{3}\rightarrow M_{\mathcal{A}}$.

\subsection{Particular results}

\noindent We present an example, originally appearing in \cite{BaMa2001} as
a fragment of the proof of the existence of a $(1,1,1,2)$ partition of two
masses by a $4$-fan on the sphere $S^{2}$.

\medskip

\noindent Let $G=\mathbb{Z}/{5}=\langle \omega \rangle $ acts on the sphere $%
S^{3}$ freely and on $\mathbb{R}^{5}$ by the cyclic shift. Again $W_{5}$ is
the subspace of $\mathbb{R}^{5}$ given by $\{(x_{1},\ldots
,x_{5})~|~x_{1}+\cdots +x_{5}=0\}$. Let $\mathcal{A}$ be a minimal $G$%
-invariant arrangement containing subspace $L\subset W_{5}$ defined by%
\begin{equation*}
x_{1}=x_{2}=x_{3}=x_{4}+x_{5}=0.
\end{equation*}%
\emph{Then there are no }$\mathbb{Z}/5$\emph{-map }$S^{3}\rightarrow
W_{5}\backslash \bigcup A=M_{\mathcal{A}}$\emph{\ and the }$\mathbb{Z}/5$%
\emph{-action on }$H^{\ast }(M_{\mathcal{A}},\Bbbk )$\emph{\ is nontrivial
for any field }$\Bbbk $\emph{.}

\noindent The Hasse diagram of the intersection poset of the arrangement $%
\mathcal{A}$ is pictured as follows.

\begin{figure}[tbh]
\centering\includegraphics[scale=0.40]{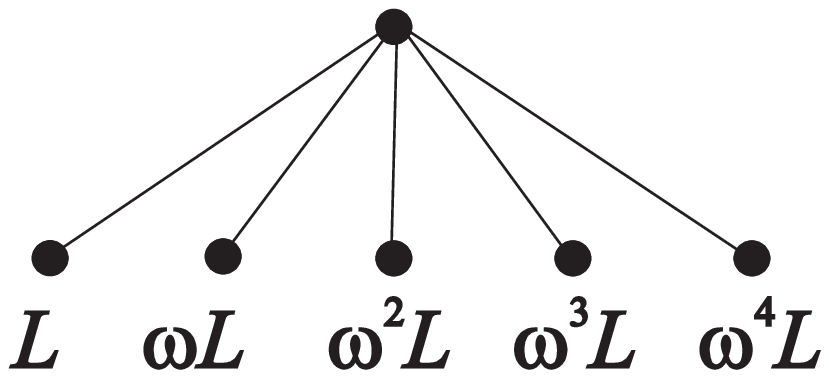}
\caption{{\protect\small Hasse diagram.}}
\label{Fig:tableu1}
\end{figure}

\noindent The cohomology of the complement $M_{\mathcal{A}}$ with
coefficients in $\mathbb{F}_{5}$ can be derived via the Goresky-MacPherson
formula (as $G$-module) as%
\begin{equation*}
H^{r}(M_{\mathcal{A}},\mathbb{F}_{5})=\left\{
\begin{array}{ll}
\mathbb{F}_{5}[\mathbb{Z}/{5}]\oplus \mathbb{F}_{5}[\mathbb{Z}/{5}%
]/(1+\omega +\cdots +\omega ^{4})\mathbb{F}_{5}[\mathbb{Z}/{5}], & r=2, \\
\mathbb{F}_{5}, & r=0, \\
0, & \text{{\small otherwise.}}%
\end{array}%
\right.
\end{equation*}%
The $E_{2}$-term of the Serre spectral sequence associated with the Borel
construction of the fibration $EG\times _{G}M_{\mathcal{A}}\rightarrow BG$
is given by%
\begin{equation*}
E_{2}^{p,q}=H^{p}(BG,\mathcal{H}^{q}(M_{\mathcal{A}},\mathbb{F}%
_{5}))=H^{p}(G,H^{q}(M_{\mathcal{A}},\mathbb{F}_{5})).
\end{equation*}%
The cohomology of the group $\mathbb{Z}/{5}$ with coefficients in the
modules
\begin{equation*}
\mathbb{F}_{5}[\mathbb{Z}/{5}],\quad \mathbb{F}_{5}[\mathbb{Z}/{5}%
]/(1+\omega +\cdots +\omega ^{4})\mathbb{F}_{5}[\mathbb{Z}/{5}]\text{\textrm{%
\ and trivial module\ }}\mathbb{F}_{5}
\end{equation*}%
is well known. 
Therefore,%
\begin{equation*}
E_{2}^{p,q}=\left\{
\begin{array}{ll}
\mathbb{F}_{5}\oplus \mathbb{F}_{5}, & q=2\text{, }p=0 \\
\mathbb{F}_{5}, & q=2\text{, }p\geq 1 \\
\mathbb{F}_{5}, & q=0\text{, }p\geq 0 \\
0, & \text{{\small otherwise}.}%
\end{array}%
\right.
\end{equation*}%
The only possibly nontrivial differential in this spectral sequence is $%
d_{3} $.

\noindent The $G$-action on $M_{\mathcal{A}}$ is free and therefore there is
a homotopy equivalence
\begin{equation*}
EG\times _{G}M_{\mathcal{A}}\simeq M_{\mathcal{A}}/G
\end{equation*}%
and consequently a group isomorphism
\begin{equation*}
H^{\ast }(EG\times _{G}M_{\mathcal{A}},\mathbb{F}_{5})\cong H^{\ast }(M_{%
\mathcal{A}}/G,\mathbb{F}_{5}).
\end{equation*}%
Since $H^{p}(M_{\mathcal{A}}/G,\mathbb{F}_{5})=0$ for $p\geq 4$, and the
fact that the Serre spectral sequence converges to $H^{\ast }(EG\times
_{G}M_{\mathcal{A}},\mathbb{F}_{5})$ the only possibly nontrivial
differential $d_{3}$ is indeed nontrivial. In particular, $d_{3}\colon
E_{3}^{0,2}\rightarrow E_{3}^{3,0}$ is nontrivial and
\begin{equation*}
E_{3}^{3,0}(EG\times _{G}M_{\mathcal{A}})\ncong E_{\infty }^{3,0}(EG\times
_{G}M_{\mathcal{A}}).
\end{equation*}

\noindent Let us assume that there is a $G$-map $S^{3}\rightarrow M_{%
\mathcal{A}}$. The induced map on the Serre spectral sequences of Borel
constructions given by%
\begin{equation*}
E_{2}^{p,0}(f)\colon E_{2}^{p,0}(EG\times _{G}M_{\mathcal{A}})\rightarrow
E_{2}^{p,0}(EG\times _{G}S^{3})
\end{equation*}%
is the identity. Since $E_{2}^{3,0}(EG\times _{G}S^{3})=E_{\infty
}^{3,0}(EG\times _{G}S^{3})$ there is an element $0\neq x\in
E_{2}^{3,0}(EG\times _{G}M_{\mathcal{A}}\times _{G}EG)$ such that
\begin{equation*}
E_{2}^{3,0}(EG\times _{G}M_{\mathcal{A}})\ni x\overset{E_{2}(f)}{\longmapsto
}x\in E_{2}^{2,0}(EG\times _{G}S^{3})
\end{equation*}%
and%
\begin{equation*}
E_{\infty }^{3,0}(EG\times _{G}M_{\mathcal{A}})\ni \left( 0=\mathrm{class}%
(x)\right) \overset{E_{\infty }(f)}{\longmapsto }\left( x\neq 0\right) \in
E_{\infty }^{2,0}(EG\times _{G}S^{3})\text{.}
\end{equation*}%
This is a \textbf{contradiction} to the assumption of the existence of a $G$%
-map $S^{3}\rightarrow M_{\mathcal{A}}$. Thus, there is no $G$-map $%
S^{3}\rightarrow M_{\mathcal{A}}$.

\section{Proof of Theorem \protect\ref{Th:APP}}

\noindent Motivated by ideas in \cite{BaMa2001} and \cite{Bl-Di}, we
consider questions of the existence of partitions and transform them to
questions of the existence of equivariant maps.

\noindent Let $k\in \mathbb{N}$ and $\mathcal{M}=\{\mu _{1},\ldots ,\mu
_{j}\}$ be a collection of measures on $S^{d-1}$. Let $\mathbf{\alpha }=(%
\tfrac{\alpha _{1}}{n},\ldots ,\tfrac{\alpha _{2k}}{n})\in \tfrac{1}{n}%
\mathbb{N}^{2k}\subset \mathbb{Q}^{2k}$ be a rationa, that is,

\begin{itemize}
\item $\alpha _{1}+\cdots +\alpha _{2k}=n$

\item for all $i\in \{1,\ldots ,k\}$, we have $\alpha _{i}=\alpha _{k+i}$.
\end{itemize}

\noindent These two conditions imply that $n$ is even.

\subsection{Configuration space}

\noindent The configuration space associated with the measure $\mu _{1}$ is
defined by%
\begin{equation*}
X_{\mu _{1},n}=\{(L,\mathcal{O}_{1},\ldots ,\mathcal{O}_{n})\in \mathcal{F}%
_{n}~|~\mu _{1}(\mathcal{O}_{1})=\cdots =\mu _{1}(\mathcal{O}_{n})=\tfrac{1}{%
n}\}.
\end{equation*}%
Any $n$-fan $(L,\mathcal{O}_{1},\ldots ,\mathcal{O}_{m})=(L;v_{1},\ldots
,v_{n})$ in the configuration space $X_{\mu _{1},n}$ is completely
determined by the vector $v_{1}$ and the orientation of the circle $%
S(L^{\bot })$. The orientation on $S(L^{\bot })$ is determined by an unit
tangent vector $v$ to the circle $S(L^{\bot })$ at a point $S(L^{\bot })\cap
\mathrm{span}\{v_{1}\}$ inside $L^{\bot }$. Therefore, we can identify%
\begin{equation*}
X_{\mu _{1},n}\cong V_{2}(\mathbb{R}^{d}).
\end{equation*}%
To see this, let $(u,w)\in V_{2}(\mathbb{R}^{d})$. Then the subspace $L$ can
be recovered as $\left( \mathrm{span}\{u,w\}\right) ^{\bot }$ and the vector
$v_{1}=u$. So far we have the subspace $L$ and the first half-hyperplane $%
F_{1}$ determined by $u$. The half-hyperplane $F_{2}$ is the one in the
direction determined by $w$ such that the measure $\mu _{1}$ of the open
sphere sector determined by $F_{1}$ and $F_{2}$ is $\tfrac{1}{n}$. The
process continues until we recover all the half-hyperplanes of a fan. There
is a natural action of the dihedral group $D_{2n}=\langle \varepsilon
,\sigma \mid \varepsilon ^{k}=\sigma ^{2}=1$,\ $\varepsilon ^{k-1}\sigma
=\sigma \varepsilon \rangle $ on $\mathcal{F}_{n}$ given by%
\begin{equation*}
\begin{array}{lll}
\varepsilon \cdot \left( L;v_{1},\ldots ,v_{n}\right) & = & \left(
L;v_{n},v_{1},\ldots ,v_{n-1}\right) , \\
\sigma \cdot \left( L;v_{1},\ldots ,v_{n}\right) & = & \left(
L;v_{n},v_{n-1},\ldots ,v_{1}\right) .%
\end{array}%
\end{equation*}

\subsection{Test Maps}

\noindent Let $W_{n}=\{(x_{1},\ldots ,x_{n})\in \mathbb{R}^{n}\mid
x_{1}+\cdots +x_{n}=0\}\subset \mathbb{R}^{n}$. A $D_{2n}$-action on $%
\mathbb{R}^{n}$ and on $W_{n}$ is given by%
\begin{equation*}
\begin{array}{lll}
\omega \cdot (x_{1},\ldots ,x_{n}) & = & (x_{n},x_{1},\ldots ,x_{n-1}), \\
\sigma \cdot (x_{1},\ldots ,x_{n}) & = & (x_{n},x_{n-1},\ldots ,x_{1}).%
\end{array}%
\end{equation*}%
The group $D_{2n}$ acts diagonally on the sum $\left( W_{n}\right) ^{\oplus
l}$.

\smallskip

\noindent \textbf{(A)} A test map $F\colon X_{\mu _{1},n}\rightarrow \left(
W_{n}\right) ^{\oplus (j-1)}$ associated with a $\mathbf{\alpha }$-partition
$k$-fan problem (Theorem \ref{Th:APP} (A)) is defined by%
\begin{equation*}
X_{\mu _{1},n}\ni (L,\mathcal{O}_{1},\ldots ,\mathcal{O}_{n})\longmapsto
\left( \mu _{i}(\mathcal{O}_{1})-\tfrac{1}{n},\ldots ,\mu _{i}(\mathcal{O}%
_{n})-\tfrac{1}{n}\right) _{i=2}^{j}\in \left( W_{n}\right) ^{\oplus (j-1)}.
\end{equation*}

\smallskip

\noindent \textbf{(B)} A test map $H:X_{\mu _{1},n}\rightarrow W_{n}\oplus
\left( W_{n}\right) ^{\oplus (j-1)}$ associated with and $\mathbf{\alpha }$%
-partition $k$-fan position arrangement is given by%
\begin{equation*}
X_{\mu _{1},n}\ni (L,\mathcal{O}_{1},\ldots ,\mathcal{O}_{n})\longmapsto
\left( \phi _{r}-\tfrac{2\pi }{n}\right) _{r=1}^{n}\times \left( \left( \mu
_{i}(\mathcal{O}_{t})-\tfrac{1}{n}\right) _{t=1}^{n}\right) _{i=2}^{j}\in
W_{n}\oplus \left( W_{n}\right) ^{\oplus (j-1)}
\end{equation*}%
where $(L,\mathcal{O}_{1},\ldots ,\mathcal{O}_{n})=(L;v_{1},\ldots ,v_{n})$
and $\phi _{r}$ denotes the angle between $v_{r}$ and $v_{r+1}$ (as
introduced in Section \ref{Sec:appResults}).

\noindent Both maps are defined in such a way that the following proposition
holds:

\begin{proposition}
The maps $F:X_{\mu _{1},n}\rightarrow \left( W_{n}\right) ^{\oplus (j-1)}$
and $H:X_{\mu _{1},n}\rightarrow W_{n}\oplus \left( W_{n}\right) ^{\oplus
(j-1)}$ are $D_{2n}$-equivariant maps.
\end{proposition}

\subsection{Test spaces}

\noindent Natural test spaces for both statements of Theorem \ref{Th:APP}
are arrangements, introduced in the following way.

\smallskip \noindent \textbf{(A)} Let $\mathcal{B}$ be the minimal $D_{2n}$%
-invariant arrangement in $\left( W_{n}\right) ^{\oplus (j-1)}$ containing
the subspace $L_{\mathcal{B}}$ given by following $k\times (j-1)$ equalities%
\begin{equation*}
\begin{array}{ll}
x_{1,i}+\cdots +x_{\alpha _{1},i} & =x_{\frac{n}{2}+1,i}+\cdots +x_{\frac{n}{%
2}+\alpha _{1},i} \\
x_{\alpha _{1}+1,i}+\cdots +x_{\alpha _{1}+\alpha _{2},i} & =x_{\frac{n}{2}%
+\alpha _{1}+1,i}+\cdots +x_{\frac{n}{2}+\alpha _{1}+\alpha _{2},i} \\
\vdots &  \\
x_{\alpha _{1}+\cdots +\alpha _{k-1}+1,i}+\cdots +x_{\frac{n}{2},i} & =x_{%
\frac{n}{2}+\alpha _{1}+\cdots +\alpha _{k-1}+1,i}+\cdots +x_{n,i}%
\end{array}%
\end{equation*}%
for all $i\in \{1,\ldots ,j-1\}$. Here $x_{a,b}$ denotes the $a$-th
coordinate in the $b$-th copy of $W_{n}$.

\smallskip \noindent \textbf{(B)} Let $\mathcal{A}$ be the minimal $D_{2n}$%
-invariant arrangement in $W_{n}\oplus \left( W_{n}\right) ^{\oplus (j-1)}$
containing the subspace $L_{\mathcal{A}}$ described by $k+k\times (j-1)$
equalities%
\begin{equation*}
\begin{array}{lll}
x_{1,1}+\cdots +x_{\frac{n}{2},1} & = & 0 \\
x_{\alpha _{1}+1,1}+\cdots +x_{\alpha _{1}+\frac{n}{2},1} & = & 0 \\
\vdots &  &  \\
x_{\alpha _{1}+\cdots +\alpha _{l-1}+1,1}+\cdots +x_{\alpha _{1}+\cdots
+\alpha _{l-1}+\frac{n}{2},1} & = & 0%
\end{array}%
\end{equation*}%
and%
\begin{equation*}
\begin{array}{ll}
x_{1,i}+\cdots +x_{\alpha _{1},i} & =x_{\frac{n}{2}+1,i}+\cdots +x_{\frac{n}{%
2}+\alpha _{1},i} \\
x_{\alpha _{1}+1,i}+\cdots +x_{\alpha _{1}+\alpha _{2},i} & =x_{\frac{n}{2}%
+\alpha _{1}+1,i}+\cdots +x_{\frac{n}{2}+\alpha _{1}+\alpha _{2},i} \\
\vdots &  \\
x_{\alpha _{1}+\cdots +\alpha _{k-1}+1,i}+\cdots +x_{\frac{n}{2},i} & =x_{%
\frac{n}{2}+\alpha _{1}+\cdots +\alpha _{k-1}+1,i}+\cdots +x_{n,i}%
\end{array}%
\end{equation*}%
for all $i\in \{2,\ldots ,j\}$. The test spaces are defined in such a way
that the following basic proposition holds.

\begin{proposition}
\label{prop:CSTM}\qquad \newline

(A) If there is no $D_{2n}$-equivariant map%
\begin{equation*}
V_{2}(\mathbb{R}^{d})\rightarrow \left( W_{n}\right) ^{\oplus
(j-1)}\backslash \bigcup \mathcal{B},
\end{equation*}%
then the statement of Theorem \ref{Th:APP} (A) is true.

(B) If there is no $D_{2n}$-equivariant map%
\begin{equation*}
V_{2}(\mathbb{R}^{d})\rightarrow W_{n}\oplus \left( W_{n}\right) ^{\oplus
(j-1)}\backslash \bigcup \mathcal{A},
\end{equation*}%
then the statement of Theorem \ref{Th:APP} (B) is true.
\end{proposition}

\subsection{Applying Theorem \protect\ref{Th:Main}}

\noindent Proposition \ref{prop:CSTM} provides a chance to apply Theorem \ref%
{Th:Main}. Unfortunately, conditions (B) and (C) for the group $D_{2n}$ with
either of arrangements $\mathcal{A}$ and $\mathcal{B}$ are not satisfied. To
overcome this difficulty we substitute the dihedral group $D_{2n}$ with its
subgroup $G=\langle \varepsilon ^{\tfrac{n}{2}}\rangle \cong \mathbb{Z}/{2}$%
. The assumptions of Theorem \ref{Th:Main} are satisfied:

\begin{itemize}
\item The Stiefel manifold $V_{2}(\mathbb{R}^{d})$ is $d-3$ connected and so
$H^{i}(V_{2}(\mathbb{R}^{d}),\mathbb{F}_{2})=0$ for $1\leq i\leq d-3$;

\item the codimension of the maximal elements of the arrangement $\mathcal{B}
$ inside $(W_{n})^{\oplus (j-1)}$ is $k(j-1)$ and the codimension of the
maximal elements in $\mathcal{A}$ inside $W_{n}\oplus \left( W_{n}\right)
^{\oplus (j-1)}$is $kj$;

\item since rations are symmetric, then $\varepsilon ^{\tfrac{n}{2}}\cdot $ $%
L_{\mathcal{A}}=L_{\mathcal{A}}$ and $\varepsilon ^{\tfrac{n}{2}}\cdot L_{%
\mathcal{B}}=L_{\mathcal{B}}$; thus, the blow-ups $\mathfrak{B}(\mathcal{A})$
and $\mathfrak{B}(\mathcal{B})$ can be constructed in such a way that $G$
acts trivially on the $\mathbb{F}_{2}$ cohomology of the complements;

\item for all $g\in D_{2n}$, we have $g\cdot L_{\mathcal{B}}\supseteq \left(
(W_{n})^{\oplus (j-1)}\right) ^{G}$ and $g\cdot L_{\mathcal{A}}\supseteq
(W_{n}\oplus \left( W_{n}\right) ^{\oplus (j-1)})^{G}$;

\item the $G$-action on the complements $\left( W_{n}\right) ^{\oplus
(j-1)}\backslash \bigcup \mathcal{B}$ and $W_{n}\oplus \left( W_{n}\right)
^{\oplus (j-1)}\backslash \bigcup \mathcal{A}$ is free.
\end{itemize}

\noindent Theorem \ref{Th:Main} implies that there are no $G$-equivariant
maps%
\begin{equation*}
V_{2}(\mathbb{R}^{d})\rightarrow \left( W_{n}\right) ^{\oplus
(j-1)}\backslash \bigcup \mathcal{B}
\end{equation*}%
and

\begin{equation*}
V_{2}(\mathbb{R}^{d})\rightarrow W_{n}\oplus \left( W_{n}\right) ^{\oplus
(j-1)}\backslash \bigcup \mathcal{A}.
\end{equation*}

\noindent and consequently there are no $D_{2n}$-equivariant maps. This
proves Theorem \ref{Th:APP}.


\begin{thebibliography}{99}
\bibitem{BaMa2001} \textsc{I.B\'{a}r\'{a}ny, J.Matou\v{s}ek,} \emph{%
Simultaneous partitions of measures by }$k$\emph{-fans}, Discrete Comp.
Geometry, 25\thinspace\ (2001), 317--334.

\bibitem{BaMa2002} \textsc{I.B\'{a}r\'{a}ny, J. Matou\v{s}ek,} \emph{%
Equipartitions of two measures by a }$4$\emph{-fan,} Discrete Comput. Geom.,
27\thinspace :\ 293-302, 2002.

\bibitem{Bl} \textsc{P. Blagojevi\'{c},} \emph{The partition of measures by }%
$3$\emph{-fans and computational obstruction theory}, arXiv: math.CO/
0402400, 2004.

\bibitem{Bl-Di} \textsc{P. Blagojevi\'{c}, A. Dimitrijevi\'{c} Blagojevi\'{c}%
}, \emph{Using Equivariant Obstruction Theory in Combinatorial Geometry},
Topology and its Applications 154 (2007) 2635-2655, %
\url{http://dx.doi.org/10.1016/j.topol.2007.04.007}.

\bibitem{Bl-Vr-Ziv} \textsc{P. Blagojevi\'{c}, S. Vre\'{c}ica, R. \v{Z}%
ivaljevi\'{c},} \emph{Computational Topology of Equivariant Maps from
Spheres to Complement of Arrangements,} arXiv:math.AT/0403161.

\bibitem{Brown} \textsc{K.S. Brown}, \emph{Cohomology of groups}, New York,
Berlin, Springer-Verlag, 1982.

\bibitem{Dold} \textsc{A. Dold,} \emph{Simple proofs of some Borsuk-Ulam
results}, Contemporary Mathematics, Vol 19, 65-69, 1983

\bibitem{Fadell-Husseini} \textsc{E.\ Fadell, S. Husseini}, \emph{An
ideal-valued cohomological index, theory with applications to Borsuk-Ulam
and Bourgin-Yang theorems}, Ergod.\ Th. and Dynam. Sys. 8$^{\ast }$(1988),
73-85.

\bibitem{Mark-Carsten} \textsc{M. de Longueville, C. A. Schultz}, \emph{The
cohomology rings of complements of subspace arrangements}, Math, Ann, 319,
625-646 (2001).

\bibitem{Mat94} \textsc{J.~Matou\v{s}ek}, \emph{Topological methods in
Combinatorics and Geometry,} Lecture notes, Prague 1994. (updated version,
February 2002, www.ms.mff.cuni.cz/matousek/lecturenotes.html).

\bibitem{McC} \textsc{J. McCleary}, \emph{A User's Guide to Spectral
Sequences}, New York, NY, Cambridge University Press, second edition, 2000.

\bibitem{SunWelker} \textsc{S. Sundaram, V. Welker,} \emph{Group actions on
arrangements of linear subspaces and applications to configuration spaces},
TAMS, Vol 349, No. 4, 1997, 1389-1420

\bibitem{Dieck87} \textsc{T. tom Dieck}, \emph{Transformation groups}, de
Gruyter Studies in Math. 8, Berlin, 1987.

\bibitem{Wall} \textsc{C.T.C. Wall,} \emph{Surgery on Compact Manifolds},
Academic Press, 1970.
\end{thebibliography}
\end{document}